\documentclass[12pt]{article}
\usepackage{amssymb}
\usepackage{amsmath}
\newtheorem{thm}{Theorem}
 
\newtheorem{la}[thm]{Lemma}
\newtheorem{cor}[thm]{Corollary} 
\newtheorem{dftemp}[thm]{Definition}
\newtheorem{extemp}[thm]{Example} 
\newtheorem{rmktemp}[thm]{Remark}
\newtheorem{convtemp}[thm]{Convention}

\newenvironment{rmk}{\begin{rmktemp}\normalfont}{\end{rmktemp}}

\newenvironment{ls}{\begin{itemize}}{\end{itemize}}
\newenvironment{lsnum}{\begin{enumerate}}{\end{enumerate}}
\newenvironment{pf}{\smallskip\noindent\emph{Proof}\quad}%
{\hfill$\square$\par\smallskip}
\newenvironment{pfof}[1]{\smallskip\noindent\emph{Proof of #1}%
\quad}{\hfill$\square$\par\smallskip}

\newcommand{\scr}[1]{\ensuremath{\mathcal {#1}}}

\newcommand{\bbb}[1]{\ensuremath{\mathbb {#1}}}

\newcommand{\emp}{\varnothing}
\renewcommand{\phi}{\varphi}

\newcommand{\notarrow}{\kern .42em\not\kern -.42em\longrightarrow}
\newcommand{\supp}{\text{supp}}
\newcommand{\pd}[2]{\ensuremath{\Pi(#1,{<}#2)}}

\date{}

\title{Ultrafilters and Partial Products\\
of Infinite Cyclic Groups}
\author{Andreas Blass%
\thanks{Partially supported by NSF grant DMS--0070723.  Part of this
  paper was written during a visit of the first author to the Centre
  de Recerca Matem\`atica in Barcelona.}\\
Mathematics Department\\
University of Michigan\\
Ann Arbor, MI 48109--1109, U.S.A.\\
ablass@umich.edu
\and
Saharon Shelah%
\thanks{Partially supported by US-Israel Binational Science
Foundation grant 2002323.  Publication number 854 of the second
author.}\\ 
Mathematics Department\\
Hebrew University\\
Jerusalem 91904, Israel\\
and\\
Mathematics Department\\
Rutgers University\\
New Brunswick, NJ 08903, U.S.A.}

\begin{document}
\maketitle

\begin{abstract}
We consider, for infinite cardinals $\kappa$ and $\alpha\leq\kappa^+$,
the group \pd\kappa\alpha\ of sequences of integers, of length
$\kappa$, with non-zero entries in fewer than $\alpha$ positions.  Our
main result tells when \pd\kappa\alpha\ can be embedded in
\pd\lambda\beta.  The proof involves some set-theoretic results, one
about families of finite sets and one about families of ultrafilters.
\end{abstract}

\section{Introduction}   \label{intro}

For an infinite cardinal $\kappa$, let $\bbb Z^\kappa$ be the direct
product of $\kappa$ copies of the additive group \bbb Z of integers.
An element of $\bbb Z^\kappa$ is thus a function\footnote{We use the
standard notational conventions of set theory, whereby a cardinal
number is an initial ordinal number and is identified with the set of
all smaller ordinals.  In particular, the cardinal of countable
infinity is identified with the set $\omega$ of natural numbers.}
$x:\kappa\to\bbb Z$, and we define its \emph{support} to be the set
$$
\supp(x)=\{\xi\in\kappa:x(\xi)\neq0\}.
$$
The partial products mentioned in the title of this paper are the
subgroups of $\bbb Z^\kappa$ of the form
$$
\pd\kappa\alpha=\{x\in\bbb Z^\kappa:|\supp(x)|<\alpha\}
$$ where $\alpha$ is an infinite cardinal no larger than the successor
cardinal $\kappa^+$ of $\kappa$.  Notice that $\pd\kappa{\kappa^+}$ is
the full product $\bbb Z^\kappa$.  At the other extreme,
\pd\kappa\omega\ is the direct sum of $\kappa$ copies of \bbb Z, i.e.,
the free abelian group generated by the $\kappa$ standard unit vectors
$e_\xi$ defined by $e_\xi(\xi)=1$ and $e_\xi(\eta)=0$ for
$\xi\neq\eta$. 

The main result in this paper gives necessary and sufficient
conditions for one partial product of \bbb Z's to be isomorphically
embeddable in another.

\begin{thm}   \label{main}
  \pd\kappa\alpha\ is isomorphic to a subgroup of \pd\lambda\beta\ if
  and only if either
  \begin{lsnum}
    \item $\kappa\leq\lambda$ and $\alpha\leq\beta$ or
    \item $\kappa\leq\lambda^{<\beta}$ and $\alpha=\omega$.
  \end{lsnum}
\end{thm}

Part of this was proved in \cite[Theorem~23 and Remark~28]{sets},
using a well-known result from set theory, the $\Delta$-system lemma.
Specifically, the results in \cite{sets} establish Theorem~\ref{main}
when $\alpha$ and $\beta$ are regular, uncountable cardinals smaller
than all measurable cardinals.  In the present paper, we complete the
proof by handling the cases of singular cardinals and cardinals above
a measurable one.  In contrast to the situation in \cite{sets}, this
will involve developing some new results in set theory, rather than
only invoking classical facts.

The set theoretic facts we need are the following two.

\begin{thm}  \label{cut}
Let $\kappa$ be an infinite cardinal and let $(F_\xi)_{\xi\in\kappa}$
be a $\kappa$-indexed family of nonempty, finite sets.  
\begin{lsnum}
  \item There exists a set $X$ such that
  \begin{equation}  \label{sing}
|\{\xi\in\kappa:|F_\xi\cap X|=1\}|=\kappa.
  \end{equation}
\item The set $X$ in \eqref{sing} can be chosen so that 
  \begin{ls}
    \item $|X|$ has cardinality 1 or $\text{cf}(\kappa)$ or $\kappa$,
    \item every subset of $X$ with the same cardinality as $X$ has the
  property in \eqref{sing}, and
  \item each element of $X$ is the unique element of $F_\xi\cap X$ for
  at least one $\xi$.
  \end{ls}
\end{lsnum}
\end{thm}

\begin{thm}  \label{ufs}
  Let $\kappa$ be an infinite cardinal and let $(\scr
  U_\xi)_{\xi\in\kappa}$ be a $\kappa$-indexed family of non-principal
  ultrafilters on $\kappa$.  Then there exists $X\subseteq\kappa$ such
  that $|X|=\kappa$ and, for each $\xi\in X$, $X\notin\scr U_\xi$.
\end{thm}

We prove these two set-theoretic theorems in Section~\ref{st}.  Then,
in Section~\ref{ab}, we apply them to prove Theorem~\ref{main}.

\section{Set Theory}   \label{st}

We thank Stevo Todor\v cevi\'c for suggesting a simplification, using
the $\Delta$-system lemma, of our original proof of Theorem~\ref{cut}.
That suggestion led us, by further simplification, to the following
proof, which doesn't need the $\Delta$-system lemma.

\begin{pfof}{Theorem~\ref{cut}}
Let $\kappa$ and $(F_\xi)_{\xi\in\kappa}$ be given, as in the
hypothesis of the theorem.  If there exists some $x$ that lies in
$F_\xi$ for $\kappa$ values of $\xi$, then $X=\{x\}$ obviously
satisfies the conclusion of the theorem.  So we assume from now on
that each $x$ lies in $F_\xi$ for fewer than $\kappa$ values of $\xi$.

\begin{la}   \label{startD}
  There exists a set $D\subseteq\kappa$ with $|D|=\kappa$ and there is
  a function assigning, to each $\xi\in D$, some $x_\xi\in F_\xi$ with
  the following property.  Whenever $x_\alpha\in F_\beta$ with
  $\alpha,\beta\in D$, then $x_\alpha=x_\beta$.
\end{la}

\begin{pf}
  We begin by simplifying a special case that would otherwise
  interfere with the main argument.  The special case is that $\kappa$
  is singular, say with cofinality $\mu$, and that there are, for
  arbitrarily large $\lambda<\kappa$, finite sets $G$ such that
  $|\{\xi<\kappa:F_\xi=G\}|=\lambda$.  That is, there are finite sets
  $G$ that are repeated nearly $\kappa$ times in the family
  $(F_\xi)_{\xi\in\kappa}$.  (Note that no set can be repeated
  $\kappa$ times, thanks to our standing assumption that no $x$ occurs
  in $F_\xi$ for $\kappa$ values of $\xi$.)  In this case, we can fix
  an increasing $\mu$-sequence of cardinals $\lambda_i$ with supremum
  $\kappa$, and we can fix finite sets $G_i$ such that each $G_i$ is
  equal to $F_\xi$ for $\lambda_i$ values of $\xi$.  Then we apply the
  argument given below to the family $(G_i)_{i\in\mu}$ instead of the
  original family $(F_\xi)_{\xi\in\kappa}$.  The result will be a set
  $D'\subseteq\mu$ of cardinality $\mu$ and a function assigning to
  each $i\in D'$ some $x'_i\in G_i$ such that, whenever $x'_i\in G_j$
  with $i,j\in D'$ then $x'_i=x'_j$.  Then we define $D$ to be the set
  of all those $\xi$ such that $F_\xi=G_i$ for some $i\in D'$, and we
  define, for each such $\xi$, $x_\xi$ to be $x'_i$, where $i\in D'$
  with $x_\xi=x'_i$.  (The defining property of $D'$ and the
  $x'_i$'s ensures that this $x'_i$ is uniquely determined for each
  $\xi$.)  It is easy to verify that $D$ and the $x_\xi$'s are as
  required by the lemma.

This completes the proof in the exceptional case, so we assume from
now on that its case hypothesis does not hold.  This implies that, for
any set $A$ of cardinality $<\kappa$, the number of $\xi$ for which
$F_\xi\subseteq A$ is also $<\kappa$.  Indeed, since $A$ has fewer
than $\kappa$ distinct finite subsets $G$, the number of $\xi$ such
that $F_\xi\subseteq A$ is the sum, over these $G$, of their
multiplicities in the sequence $(F_\xi)_{\xi\in\kappa}$.  These
multiplicities are all $<\kappa$, so the only way their sum, over the
fewer than $\kappa$ $G$'s, can be $\kappa$ is for the hypothesis of
the exceptional case to hold.

We are now ready to start building the required $D$ and the required
function $\xi\mapsto x_\xi$ inductively.  We begin with $D$ empty, and
we enlarge it step by step, stopping when its cardinality reaches
$\kappa$.  At each step, we shall choose a suitable $x$ and add to
$D$ all those $\xi$ such that $x\in F_\xi$; for each of these $\xi$,
we shall set $x_\xi=x$.  In order for this definition to be consistent
and to satisfy the requirements of the lemma, our choice of $x$ is
subject to several constraints:
\begin{ls}
  \item $x$ is not in $F_\eta$ for any $\eta$ previously put into $D$.
  \item No $F_\xi$ contains both $x$ and any $x_\eta$ for $\eta$
  previously put into $D$.
  \item $x$ is in $F_\xi$ for some $\xi$.
\end{ls}
The first of these constraints ensures that the requirement in the
lemma is satisfied when $\alpha$ is one of the $\xi$'s being added at
the current step and $\beta$ was put into $D$ earlier.  The second
ensures the requirement of the lemma when $\beta$ is one of the
$\xi$'s being added at the current step and $\alpha$ was put into $D$
earlier.  (In both cases, we ensure that $x_\alpha\notin F_\beta$.)
The requirement of the lemma will also hold when both $\alpha$ and
$\beta$ are among the currently added $\xi$'s, because then
$x_\alpha=x_\beta=x$.  The third constraint merely ensures that $D$
acquires at least one new element per step; any $\xi$ as in the third
constraint is put into $D$, and it wasn't previously in $D$ because of
the first constraint.

To complete the proof of the lemma, we must show that, as long as
$|D|<\kappa$, we can find an $x$ satisfying all the constraints.  

In fact, the second constraint is redundant.  If $F_\xi$ and $\eta$
violated it, then $\xi$ would have been put into $D$ already at the
same step where $\eta$ was added, because we always add all $F$'s that
contain the currently chosen $x$.  Thus, the first constraint would be
violated with $\xi$ in the role of $\eta$.  So we need only show that,
when $|D|<\kappa$, we can choose $x$ so as to satisfy the first and
third constraints.  The union of the $F_\eta$'s for $\eta$ previously
put into $D$ is a set $A$ of cardinality $<\kappa$, because
$|D|<\kappa$ and the $F_\eta$'s are finite.  We saw above that such an
$A$ cannot include $F_\xi$ for $\kappa$ values of $\xi$.  So we can
choose a $\xi<\kappa$ with $F_\xi\not\subseteq A$ and we can choose
$x\in F_\xi-A$.  This $x$ clearly satisfies the first and third
constraints, so the proof of the lemma is complete.
\end{pf}

Fix $D$ and $\xi\mapsto x_\xi$ as in the lemma. We next normalize the
$D$ a bit as follows.  Let $\sim$ be the equivalence relation on $D$
defined by
$$
\xi\sim\eta\iff x_\xi=x_\eta.
$$
We shall arrange that one of the following three alternatives holds.
\begin{lsnum}
  \item $D$ is a single equivalence class, i.e., all the $x_\xi$ are
  equal. 
  \item Each equivalence class is a singleton, i.e., all the $x_\xi$
  are distinct.
  \item $\kappa$ is singular, the number of equivalence classes is
  $\mu=\text{cf}(\kappa)$, and their sizes form a cofinal subset of
  $\kappa$ of order-type $\mu$.
\end{lsnum}
We can arrange this simply by shrinking $D$ (while keeping its
cardinality equal to $\kappa$ of course).  If there is an equivalence
class of size $\kappa$, then replacing $D$ by that equivalence class
attains alternative (1).  If there are $\kappa$ equivalence classes,
then replacing $D$ by a selector attains alternative (2). So we may
assume that there are $<\kappa$ equivalence classes, each of size
$<\kappa$.  Thus, $\kappa$ is singular; let $\mu$ be its cofinality.
The sizes of the equivalence classes must be unbounded below $\kappa$,
for otherwise their union would be smaller than $\kappa$ (being at
most the bound times $\mu$).  So we can choose a $\mu$-sequence of
equivalence classes of increasing cardinalities approaching $\kappa$.
Replacing $D$ by the union of these classes attains alternative (3).

Finally, we let $X=\{x_\xi:\xi\in D\}$ and we check that it has the
properties required in the theorem.  If $\beta\in D$, then $x_\beta\in
F_\beta$ and, by the requirement in the lemma, no $x_\alpha\neq
x_\beta$ can be in $F_\beta$.  So $|F_\beta\cap X|=1$ for all
$\beta\in D$.  Since $|D|=\kappa$, part~1 of the theorem is satisfied.
The cardinality of $X$ is the number of equivalence classes with
respect to $\sim$ in $D$, and our normalization of $D$ ensures that
this is $1$ or $\kappa$ or $\text{cf}(\kappa)$.  The normalization
also ensures that any subset of $X$ of the same cardinality as $X$
arises from a subset of $D$ that shares the properties we obtained for
$D$.  So any such subset also works in part~1 of the theorem.
Finally, each element $x\in X$ is of the form $x_\xi$ for some $\xi\in
D$ and therefore is, thanks to the requirement on $D$ in the lemma,
the unique element of $F_\xi\cap X$.
\end{pfof}

\begin{pfof}{Theorem~\ref{ufs}}
Let $\kappa$ and $(\scr U_\xi)_{\xi\in\kappa}$ be as in the hypothesis
of the theorem.  Partition $\kappa$ into $\kappa$ sets $A_\mu$ (with
$\mu\in\kappa$), each of cardinality $\kappa$.  If one of these
$A_\mu$ can serve as $X$ in the conclusion of the theorem, then
nothing more needs to be done.  So assume that this is not the case,
i.e., assume that, for each $\mu$, there is some $\xi(\mu)\in A_\mu$
such that $A_\mu\in\scr U_{\xi(\mu)}$.  Being non-principal, $\scr
U_{\xi(\mu)}$ also contains $A_\mu-\{\xi(\mu)\}$.  

Let $X=\{\xi(\mu):\mu\in\kappa\}$.  For each element of $X$, say
$\xi(\mu)$, we have seen that $\scr U_{\xi(\mu)}$ contains a set
disjoint from $X$, namely $A_\mu-\{\xi(\mu)\}$.  Therefore
$X\notin\scr U_{\xi(\mu)}$, and the proof is complete.
\end{pfof}

\section{Proof of Theorem~\ref{main}}   \label{ab}

We begin by showing that, if one of the cardinality conditions 1
and 2 in Theorem~1 is satisfied, then we can embed \pd\kappa\alpha\ in
\pd\lambda\beta.  

If $\kappa\leq\lambda$, then we can embed $\bbb Z^\kappa$ into $\bbb
Z^\lambda$ by extending any $\kappa$-sequence $x\in\bbb Z^\kappa$ by
zeros to have length $\lambda$.  This does not alter the support, so
it embeds \pd\kappa\alpha\ into \pd\lambda\beta\ (as a pure subgroup)
for any $\beta\geq\alpha$.

This completes the proof if condition 1 in the theorem is satisfied.
If condition 2 is satisfied, then, since $\alpha=\omega$, the group
\pd\kappa\alpha\ is a free abelian group of rank
$\kappa\leq\lambda^{<\beta}$.  Since \pd\lambda\beta\ has cardinality
$\lambda^{<\beta}$, its rank is also $\lambda^{<\beta}$.  (The only
way for a torsion-free abelian group to have rank different from its
cardinality is to have finite rank, which is clearly not the case for
\pd\lambda\beta.)  So it has a free subgroup of rank
$\lambda^{<\beta}$, and we have the required embedding.

\begin{rmk}
  N\"obeling proved in \cite{nob} that the subgroup of $\bbb
  Z^\lambda$ consisting of the bounded functions is a free abelian
  group.  Intersecting it with \pd\lambda\beta, we get a pure free
  subgroup of \pd\lambda\beta\ of rank $\lambda^{<\beta}$.  Thus, under
  condition~2 of the theorem, we get an embedding of \pd\kappa\alpha\
  into \pd\lambda\beta\ as a pure subgroup.  Therefore, Theorem~1
  would remain correct if we replaced ``subgroup'' with ``pure
  subgroup.''
\end{rmk}

We now turn to the more difficult half of Theorem~\ref{main}, assuming
the existence of the embedding of groups and deducing one of the
cardinality conditions.  Since \pd\lambda\beta\ has cardinality
$\lambda^{<\beta}$ and \pd\kappa\alpha\ has cardinality at least
$\kappa$, the existence of an embedding of the latter into the former
obviously implies that $\kappa\leq\lambda^{<\beta}$.  So if
$\alpha=\omega$ then we have condition~2 of the theorem.
Therefore, we assume from now on that $\alpha$ is uncountable; our
goal is to deduce condition~1.

For this purpose, we need to assemble some information about the given
embedding $j:\pd\kappa\alpha\to\pd\lambda\beta$.  The embedding is, of
course, determined by its $\lambda$ components, i.e., its compositions
with the $\lambda$ projection functions $p_\nu:\pd\lambda\beta\to\bbb
Z$.  (Here and in all that follows, the variable $\nu$ is used for
elements of $\lambda$.)  We write $j_\nu$ for $p_\nu\circ
j:\pd\kappa\alpha\to\bbb Z$.  Thus, for any $x\in\pd\kappa\alpha$,
$j_\nu(x)$ is the $\nu^{\text{th}}$ component of the
$\lambda$-sequence $j(x)$.

The structure of homomorphisms, like $j_\nu$, from \pd\kappa\alpha\ to
\bbb Z can be determined, thanks to the following theorem of Balcerzyk
\cite{bal}.  (This theorem extends earlier results of Specker
\cite{specker} for $\kappa=\omega$ and \L o\'s (see
\cite[Theorem~94.4]{fuchs2}) for $\kappa$ smaller than all measurable
cardinals; it was in turn extended by Eda \cite{eda} to allow
arbitrary slender groups in place of \bbb Z.)  To state it, we need
one piece of notation.

If \scr U is a countably complete ultrafilter on a set $A$ and if $x$
is any function from $A$ to a countable set (such as \bbb Z), then
$x$ is constant on some set in \scr U, and we denote that constant
value by $\scr U\text{-}\!\lim x$.

\begin{thm}[Balcerzyk]  \label{ke}
  Let $A$ be any set and let $h:\bbb Z^A\to\bbb Z$ be a homomorphism.
  Then there exist finitely many countably complete ultrafilters $\scr
  U_i$ on $A$ and there exist integers $c_i$ (indexed by the
  same finitely many $i$'s) such that, for all $x\in\bbb Z^A$,
$$
h(x)=\sum_ic_i\cdot\scr U_i\text{-}\!\lim x.
$$
\end{thm}

We shall refer to the sum in this theorem as the \emph{Balcerzyk
formula} for $h$.  Whenever it is convenient, we shall assume that, in a
Balcerzyk formula, all the $\scr U_i$ are distinct and all the $c_i$ are
non-zero.  This can be arranged simply by combining any terms that
involve the same ultrafilter and omitting any terms with zero
coefficients.
  
The theorem easily implies that the group of homomorphisms from $\bbb
Z^A$ to \bbb Z is freely generated by the homomorphisms $\scr
U\text{-}\!\lim$ for countably complete ultrafilters \scr U on $A$.

Notice that among the countably complete ultrafilters are the
principal ultrafilters, and that the homomorphism $\scr U\text{-}\!\lim$
associated to the principal ultrafilter \scr U at some $a\in A$ is
simply the projection $p_a:\bbb Z^A\to\bbb Z:x\to x(a)$.  If $|A|$ is
smaller than all measurable cardinals, then the principal
ultrafilters are the only countably complete ultrafilters on $A$, so
homomorphisms from $\bbb Z^A$ to \bbb Z are simply finite linear
combinations of projections.

\begin{cor}   \label{finite1}
  If $h:\bbb Z^A\to\bbb Z$ is a homomorphism, then there are only
  finitely many $a\in A$ such that the standard unit vector $e_a$ is
  mapped to a non-zero value by $h$.
\end{cor}

\begin{pf}
  For $h(e_a)$ to be non-zero, one of the $\scr U_i$ in the theorem
  must be the principal ultrafilter at $a$.
\end{pf}

We wish to apply this information to the homomorphisms $j_\nu$, whose
domain is only \pd\kappa\alpha, not all of $\bbb Z^\kappa$.
Fortunately, the preceding corollary carries over to the desired
context, thanks to our assumption above that $\alpha$ is uncountable.

\begin{cor}   \label{finite2}
  For each $\nu\in\lambda$, there are only finitely many
  $\xi\in\kappa$ such that $j_\nu(e_\xi)\neq0$.
\end{cor}

\begin{pf}
  Suppose not.  Then there is a countably infinite set
  $A\subseteq\kappa$ such that, for each $\xi\in A$,
  $j_\nu(e_\xi)\neq0$.  View $\bbb Z^A$ as a subgroup of $\bbb
  Z^\kappa$, simply by extending functions by 0 on $\kappa-A$.  Since
  $\alpha$ is uncountable, we have made $\bbb Z^A$ a subgroup of
  \pd\kappa\alpha, the domain of $j_\nu$.  So we can apply
  Corollary~\ref{finite1} to (the restriction to $\bbb Z^A$ of)
  $j_\nu$ and conclude that $j_\nu(e_\xi)\neq0$ for only finitely many
  $\xi\in A$.  This contradicts our choice of $A$.
\end{pf}

For each $\nu\in\lambda$, let
$$
F_\nu=\{\xi\in\kappa:j_\nu(e_\xi)\neq0\}.
$$
So each $F_\nu$ is finite.  On the other hand, since $j$ is an
embedding, we have, for each $\xi\in\kappa$, that $j(e_\xi)\neq0$ and
therefore $\xi\in F_\nu$ for at least one $\nu\in\lambda$.  Thus,
$\kappa$ is the union of the $\lambda$ finite sets $F_\nu$, which
implies that $\kappa\leq\lambda$.  This proves the first part of
condition~1 of the theorem.

Before turning to the second part, we note, since we shall need it
later, that the preceding argument shows not only that
$\kappa\leq\lambda$ but that
$$
\kappa\leq|\{\nu\in\lambda:F_\nu\neq\emp\}|.
$$

To complete the proof of condition 1 of the theorem, it remains to
show that $\alpha\leq\beta$.  Suppose, toward a contradiction, that
$\beta<\alpha$.  So $\beta^+\leq\alpha\leq\kappa^+$ and therefore
$\beta\leq\kappa$.  Therefore (by the first part of this proof), $\bbb
Z^\beta=\pd{\beta}{\beta^+}$ embeds in \pd\kappa\alpha, which in turn
embeds in \pd\lambda\beta.  So instead of dealing with an embedding
$\pd\kappa\alpha\to\pd\lambda\beta$, we can deal with an embedding
$j:\bbb Z^\beta\to\pd\lambda\beta$.  In other words, we can assume,
without loss of generality, that $\kappa=\beta$ and $\alpha=\beta^+$.

We record for future reference that we have already reached a
contradiction if $\beta=\omega$, for then \pd\lambda\beta\ is the free
abelian group on $\lambda$ generators while, by a theorem of Specker
\cite{specker}, $\bbb Z^\beta$ is not free.  So the latter cannot be
embedded into the former.  Thus, we may assume, for the rest of this
proof, that $\beta$ is uncountable.

As before, we write $j_\nu$ for the homomorphism $\bbb Z^\beta\to\bbb
Z$ given by the $\nu^{\text{th}}$ component of $j$, for each
$\nu\in\lambda$.  Also as before, we write $F_\nu$ for the set of
$\xi\in\beta$ such that $j_\nu(e_\xi)\neq0$.  It will be useful to
write the Balcerzyk formula for $j_\nu$ with the principal and
non-principal ultrafilters separated.  Note that the principal
ultrafilters that occur here are concentrated at the points of
$F_\nu$.  Thus, we have
\begin{equation}  \label{jfmla}
  j_\nu(x)=\sum_{\xi\in F_\nu}a^\nu_\xi\cdot x(\xi)+
\sum_{\scr U\in\bbb U_\nu}b^\nu_{\scr U}\cdot\scr U\text{-}\!\lim x
\end{equation}
where $\bbb U_\nu$ is a finite set of non-principal, countably
complete ultrafilters on $\beta$.  As before, we assume, without loss
of generality, that all the $a$ and $b$ coefficients are non-zero.

We recall that we showed, in the proof of $\kappa\leq\lambda$, that
$F_\nu\neq\emp$ for at least $\beta$ values of $\nu$ (since the
$\kappa$ of that proof is now equal to $\beta$).  So we can apply
Theorem~\ref{cut} to find an $X\subseteq\beta$ with the following
properties.
\begin{lsnum}
  \item There are $\beta$ values of $\nu$, which we call the
  \emph{special} values, such that $X\cap F_\nu$ is a singleton.
  \item $|X|$ is one of 1, $\text{cf}(\beta)$, and $\beta$.
    \item Every subset of $X$ of the same cardinality as $X$ shares
    with $X$ the property in item 1 above.
    \item Each $\xi\in X$ is, for at least one $\nu$, the unique
    element of $X\cap F_\nu$.
\end{lsnum}
It will be useful to select, for each $\xi\in X$, one $\nu$ as in item
4 and to call it $\nu(\xi)$.  Notice that $\nu(\xi)$ is always special
(as defined in item 1).  

In the course of the proof, we will occasionally replace $X$ by a
subset of the same cardinality, relying on property 3 of $X$ to ensure
that all the properties listed for $X$ remain correct for the new
$X$.  To avoid an excess of subscripts, we will not give these $X$'s
different names.  Rather, at each stage of the proof, $X$ will refer
to the current set, which may be a proper subset of the original $X$
introduced above.

The basic idea of the proof is quite simple, so we present it first
and afterward indicate how to handle all the issues that arise in its
application.  

Consider any $x\in\bbb Z^\beta$ whose support is exactly $X$.  Then
for each special $\nu$ the first sum in \eqref{jfmla} reduces to a
single term, because exactly one $\xi\in F_\nu$ has $x(\xi)\neq0$.  So
this formula reads
\begin{equation}  \label{jfmla1}  
j_\nu(x)=a^\nu_\xi\cdot x(\xi)+
\sum_{\scr U\in\bbb U_\nu}b^\nu_{\scr U}\cdot\scr U\text{-}\!\lim x
\end{equation}
where $\xi$ is the unique element of $X\cap F_\nu$.  If we knew that
none of the ultrafilters $\scr U\in\bbb U_\nu$ contain $X$, then all
the corresponding limits $\scr U\text{-}\!\lim x$ would vanish, since
\scr U contains a set (namely the complement of $X$) on which $x$ is
identically 0.  In this case, we would have 
$$
j_\nu(x)=a^\nu_\xi\cdot x(\xi)\neq0.
$$
If this happened for $\beta$ distinct values of $\nu$, then all these
values would be in the support of $j(x)$, contradicting the fact that
$j(x)\in\pd\lambda\beta$.  

This is the basic idea; the rest of the proof is concerned with the
obvious difficulty that we do not immediately have $\beta$ values of
$\nu$ for which the ultrafilters $\scr U\in\bbb U_\nu$ do not contain
$X$.  

Of course, this difficulty cannot arise if $|X|=1$, as the
ultrafilters in question are non-principal.  So the proof is complete
if there is some $\xi$ that lies in $\beta$ of the sets $F_\nu$, for
then $\{\xi\}$ could serve as $X$.  From now on, we assume that there
is no such $\xi$.

More generally, the difficulty cannot arise, and so the proof is
complete, if $|X|$ is smaller than all measurable cardinals, because
then there are no non-principal, countably complete ultrafilters to
contribute to the second sum in \eqref{jfmla}.  So we may assume that
there is at least one measurable cardinal $\leq|X|$.

There remain the cases that $|X|=\beta$ and that
$|X|=\text{cf}(\beta)<\beta$.  It turns out to be necessary to
subdivide the former case according to whether
$\text{cf}(\beta)=\omega$ or not.  We handle the three resulting cases
in turn.

\smallskip\noindent\emph{Case 1:}\enspace
$|X|=\beta$ and $\text{cf}(\beta)>\omega$.

Recall that we chose, for each $\xi\in X$, some $\nu(\xi)$ such that
$X\cap F_{\nu(\xi)}=\{\xi\}$.  Thus, equation \eqref{jfmla1} holds
when we put $\nu(\xi)$ in place of $\nu$.  

There are only countably many possible values for $|\bbb
U_{\nu(\xi)}|$ because these cardinals are finite.  Since $|X|$ has,
by the case hypothesis, uncountable cofinality, $X$ must have a
subset, of the same cardinality $\beta$, such that $|\bbb
U_{\nu(\xi)}|$ has the same value, say $l$, for all $\xi$ in this
subset.  Replace $X$ with this subset; as remarked above, we do not,
with this replacement, lose any of the properties of $X$ listed above.
Now we can, for each $\xi$ in (the new) $X$, enumerate $\bbb
U_{\nu(\xi)}$ as $\{\scr U_k(\xi):k<l\}$.

Next, apply Theorem~\ref{ufs} $l$ times in succession, starting with
the current $X$.  At step $k$ (where $0\leq k<l$), replace the then
current $X$ with a subset, still of cardinality $\beta$, such that,
for each $\xi$ in (the new) $X$, $\scr U_k(\xi)$ does not contain $X$.
Thus, for the final $X$, after these $l$ shrinkings, we have that, for
all $\xi\in X$, and all $\scr U\in\bbb U_{\nu(\xi)}$, $X\notin\scr
U$.  This is exactly what we need in order to apply the basic idea,
explained above, to all the $\nu$'s of the form $\nu(\xi)$ for $\xi\in
X$.  Since the function $\xi\mapsto\nu(\xi)$ is obviously one-to-one,
there are $\beta$ of these $\nu$'s, and so we have the required
contradiction. 

Notice that the case hypothesis that $\beta$ has uncountable
cofinality was used in order to get a single cardinal $l$ for $|\bbb
U_{\nu(\xi)}|$, independent of $\xi$, which was used in turn to fix
the number of subsequent shrinkings of $X$.  Without a fixed $l$,
there would be no guarantee of a final $X$ to which the basic idea can
be applied.  This is why the following case must be treated
separately.  It is the only case where the actual values of $x$,
not just its support, will matter.

\smallskip\noindent\emph{Case 2:}\enspace
$|X|=\beta$ and $\text{cf}(\beta)=\omega$.

Recall that we have already obtained a contradiction when
$\beta=\omega$, so in the present case $\beta$ is a singular cardinal.
Fix an increasing $\omega$-sequence $(\beta_n)_{n\in\omega}$ of
uncountable regular cardinals with supremum $\beta$.  Partition $X$
into countably many sets $X_n$ with $|X_n|=\beta_n$.  As in the proof
of Case 1, we can shrink each $X_n$, without decreasing its
cardinality, so that:
\begin{ls}
\item The cardinality of $\bbb U_{\nu(\xi)}$ depends only on $n$, not
  on the choice of $\xi\in X_n$; call this cardinality $l(n)$.
  \item For all $\xi\in X_n$, no ultrafilter in $\bbb U_{\nu(\xi)}$
contains $X_n$.
\end{ls}
Here and below, when we shrink the $X_n$'s, it is to be understood
that $X$ is also shrunk, to the union of the new $X_n$'s.  As long as
the cardinality of each $X_n$ remains $\beta_n$, the cardinality of
$X$ remains $\beta$.

As before, we use the notation $\{\scr U_k(\xi):k<l(n)\}$ for an
enumeration of $\bbb U_{\nu(\xi)}$ when $\xi\in X_n$.  

Notice that each $\scr U_k(\xi)$, being countably complete, must
concentrate on one $X_m$ or on the complement of $X$.  Shrinking each
$X_n$ again without reducing its cardinality, we arrange that for each
fixed $n$ and each fixed $k<l(n)$, as $\xi$ varies over $X_n$, all the
ultrafilters $\scr U_k(\xi)$ that contain $X$ also contain the same
$X_m$.  We write $m(n,k)$ for this $m$.  (If none of these $\scr
U_k(\xi)$ contain $X$, define $m(n,k)\in\omega-\{n\}$ arbitrarily.)
Also, define $S(n)=\{m(n,k):k<l(n)\}$.  Thus, when $\xi\in X_n$, every
ultrafilter in $\bbb U_{\nu(\xi)}$ that contains $X$ contains $X_m$
for some $m\in S(n)$.  Note that our previous shrinking of the $X_n$'s
ensures that $n\notin S(n)$.

(A technical comment: When we shrink $X$ by shrinking all the $X_n$'s,
the property of an ultrafilter that ``$X_m\in\scr U$'' may be lost,
since $X_m$ may shrink to a set not in \scr U.  But, if this happens,
then $X$ also shrinks to a set not in \scr U.  Thus, the property ``if
$X\in\scr U$ then $X_m\in\scr U$'' persists under such shrinking.
This fact was tacitly used in the shrinking process of the preceding
paragraph.  It ensures that we can base our decision of how to shrink
the $X_n$'s on our knowledge of which $X_m$'s are in which
ultrafilters, without worrying that the shrinking will alter that
knowledge in a way that requires us to revise the shrinking.)

Obtain an infinite subset $Y$ of $\omega$ by choosing its elements
inductively, in increasing order, so that whenever $n<n'$ are in $Y$
then $n'\notin S(n)$.  This is trivial to do, since each $S(n)$ is
finite.  Shrink $X_n$ to $\emp$ for all $n\notin Y$, but leave $X_n$
unchanged for $n\in Y$.  Unlike previous shrinkings, this obviously
does not maintain $|X_n|=\beta_n$ in general but only for $n\in Y$.
That is, however, sufficient to maintain $|X|=\beta$, since $Y$ is
cofinal in $\omega$ and so the $\beta_n$ for $n\in Y$ have supremum
$\beta$.  As a result of this last shrinking, we have that, for each
$n\in Y$ and each $\xi\in X_n$, each of the ultrafilters $\scr
U_k(\xi)\in\bbb U_{\nu(\xi)}$ that contains $X$ also contains $X_m$
with $m=m(n,k)<n$.

Shrinking the surviving $X_n$'s further, without reducing their
cardinalities, we can arrange that in formulas \eqref{jfmla} and
\eqref{jfmla1} the coefficient $b^{\nu(\xi)}_{\scr U_k(\xi)}$ depends
only on $n$ and $k$, not on the choice of $\xi\in X_n$. We call this
coefficient $b(n,k)$.

We shall now define a certain $x\in\bbb Z^\beta$ with support (the
current) $X$.  It will be constant on each $X_n$ with a value $z_n$ to
be specified, by induction on $n$.  (Here $n$ ranges over $Y$, since
$X_n=\emp$ for $n\notin Y$.)  Suppose that integers $z_m$
have already been defined for all $m<n$.  Then for $\xi\in X_n$ the
sum in formula \eqref{jfmla1} for $\nu=\nu(\xi)$ is
$$
\sum_{\scr U\in\bbb U_{\nu(\xi)}}b^{\nu(\xi)}_{\scr U}\cdot\scr
U\text{-}\!\lim x=
\sum_{k<l(n)}b(n,k)\cdot\scr U_k(\xi)\text{-}\!\lim x=
\sum_{k<l(n)}b(n,k)\cdot(z_{m(n,k)}|0).
$$ 
Here $(z|0)$ means $z$ or 0, according to whether $\scr U_k(\xi)$
contains $X$ (and therefore $X_{m(n,k)}$) or not.  So this sum has
only finitely many (at most $2^{l(n)}$) possible values.  Choose $z_n$
to be an integer greater than the absolute values of these finitely
many possible sums.  This choice ensures that, in formula
\eqref{jfmla1} for $\nu=\nu(\xi)$ and $\xi\in X_n$, the first term
$a^{\nu(\xi)}_\xi x(\xi)$ exceeds in absolute value the sum over
non-principal ultrafilters.  Therefore, $j_{\nu(\xi)}(x)\neq0$.

But this happens for all $\xi\in X$, so $\supp(j(x))$ has cardinality
$\beta$, contrary to the fact that $j(x)\in\pd\lambda\beta$.  This
contradiction completes the proof for Case 2. 

\smallskip\noindent\emph{Case 3:}\enspace
$|X|=\text{cf}(\beta)<\beta$.

We already observed that the basic idea suffices to complete the proof
if $|X|$ is smaller than all measurable cardinals.  So in the present
situation, we may assume that $\text{cf}(\beta)$ is greater than or
equal to the first measurable cardinal; in particular it is
uncountable.  

Let $\mu=\text{cf}(\beta)$ and let $(\beta_i)_{i\in\mu}$ be an
increasing $\mu$-sequence of regular, uncountable cardinals with
supremum $\beta$.  

For each $i\in\mu$, there is some $\xi_i\in X$ such that
$$
|\{\nu:X\cap F_\nu=\{\xi_i\}\}|\geq\beta_i.
$$
Indeed, if there were no such $\xi_i$, then $\{\nu:|X\cap F_\nu|=1\}$
would be the union of $|X|=\mu$ sets each of size $<\beta_i$, so it
would have cardinality at most $\mu\cdot\beta_i<\beta$, contrary to
our original choice of $X$.

Fix such a $\xi_i$ for each $i\in\mu$.  Note that $|\{\nu:X\cap
F_\nu=\{\xi_i\}\}|$, though at least $\beta_i$ by definition, cannot
be as large as $\beta$, as we remarked when we disposed of the case
$|X|=1$ long ago.  So, although the same element can serve as $\xi_i$
for several $i$'s, it cannot do so for cofinally many $i\in\mu$.  So
there are $\mu$ distinct $\xi_i$'s.  Passing to a subsequence and
re-indexing, we henceforth assume that all the $\xi_i$ are distinct.

Next, fix for each $i\in\mu$ a set $N_i\subseteq\lambda$ of size
$\beta_i$ such that all elements $\nu$ of $N_i$ have $X\cap
F_\nu=\{\xi_i\}$.  Note that the sets $N_i$ are pairwise disjoint.

Shrink $X$ to $\{\xi_i:i\in\mu\}$.  This still has cardinality $\mu$ and
thus has all the properties originally assumed for $X$.  

For each $i$, shrink $N_i$, without reducing its cardinality
$\beta_i$, so that as $\nu$ varies over $N_i$, the cardinality of
$\bbb U_\nu$ remains constant, say $l(i)$.  This shrinking is possible
because $\text{cf}(\beta_i)>\omega$.  Since $\mu$ is uncountable and
regular, we can shrink $X$, without reducing its cardinality, so that
$l(i)$ is the same number $l$ for all $\xi_i\in X$.  Again, re-index
$X$ as $\{\xi_i:i\in\mu\}$ and re-index the $\beta_i$ and $N_i$
correspondingly.  So we can, for each $\nu\in\bigcup_iN_i$, enumerate
$\bbb U_\nu$ as $\{\scr U_k(\nu):k<l\}$.

For each $i$, choose a uniform ultrafilter $\scr V_i$ on $N_i$, and
define an ultrafilter $\scr W_i$ as the limit with respect to $\scr
V_i$ of the ultrafilters $\scr U_0(\nu)$.  That is,
$$
A\in\scr W_i\iff\{\nu:A\in\scr U_0(\nu)\}\in\scr V_i.
$$ 
It is well known and easy to check that this $\scr W_i$ is indeed an
ultrafilter.  Applying Theorem \ref{ufs}, we obtain $Y\subseteq X$ of
cardinality $\mu$, such that for each $\xi_i\in Y$, $Y\notin\scr W_i$.
This means, by definition of $\scr W_i$, that we can shrink $N_i$ to a
set in $\scr V_i$, hence still of size $\beta_i$ as $\scr V_i$ is
uniform, so that for all $\nu$ in the new $N_i$, $\scr U_0(\nu)$
doesn't contain $Y$.  Shrink $X$ to $Y$ and reindex as before. We have
achieved that, for all $i$ and all $\nu\in N_i$, $X\notin\scr
U_0(\nu)$.

Repeat the process with the subscript 0 of \scr U replaced in turn by
$1,2,\dots,l-1$.  At the end, we have $X$ and $N_i$'s such that, for all
$\xi_i\in X$, all $\nu\in N_i$, and all $\scr U\in\bbb U_\nu$,
$X\notin\scr U$.  

This means that, in formula \eqref{jfmla1} for $\xi=\xi_i\in X$ and
$\nu\in N_i$, if $x$ has support $X$, then the sum over non-principal
ultrafilters vanishes and we reach a contradiction as in the basic
idea.

\end{document}